\documentclass[11pt]{article}
\usepackage{indentfirst}

\usepackage{float}
\usepackage{booktabs}
\usepackage{makecell}
\usepackage{epsfig}
\usepackage{url}
\usepackage{ulem}
\usepackage{amsmath}
\usepackage{amsfonts}
\usepackage{color}
\usepackage{amssymb}
\usepackage{graphicx,amsfonts,mathrsfs}
\usepackage{epstopdf}
\usepackage{enumerate}
\usepackage[all]{xy}
\usepackage{amsthm}
\usepackage{array}
\usepackage[noadjust]{cite}

\allowdisplaybreaks

\begin{document}
\begin{sloppypar}

\newtheorem{theorem}{Theorem}[section]
\newtheorem{remark}[theorem]{Remark}
\newtheorem{problem}[theorem]{Problem}
\newtheorem{corollary}[theorem]{Corollary}
\newtheorem{definition}[theorem]{Definition}
\newtheorem{conjecture}[theorem]{Conjecture}
\newtheorem{question}[theorem]{Question}

\newtheorem{lemma}[theorem]{Lemma}
\newtheorem{proposition}[theorem]{Proposition}
\newtheorem{quest}[theorem]{Question}
\newtheorem{example}[theorem]{Example}
\newtheorem{observation}[theorem]{Observation}
\newcommand{\pp}{{\it p.}}
\newcommand{\de}{\em}

\title{  {On subgroup regular sets in generalized Cayley graphs}\thanks{
Weijun Liu is partly supported by the National Natural Science Foundation of China (No. 12471022),
Peng-Li Zhang is partly supported by the National Natural Science Foundation of China (No. 12501493).
 E-mail addresses: liaoqianfen@163.com(Q.F. Liao),  wjliu6210@126.com(W.J. Liu), zpengli0626@163.com(P.-L. Zhang, corresponding author).}}

\author{Qianfen Liao$^1$, Weijun Liu$^{2,3}$,  Peng-Li Zhang$^4$\\
{\small $^{1}$ Department of Mathematics, Guangdong University of Education, } \\
{\small  Guangzhou,  Guangdong, 510303, P.R. China. }\\
{\small $^{2}$ College of General Education, Guangdong University of Science and Technology,} \\
{\small Dongguan, Guangdong, 523083, P.R. China}\\
{\small $^{3}$ School of Mathematics and Statistics, Central South University} \\
{\small New Campus, Changsha, Hunan, 410083, P.R. China. }\\
{\small $^{4}$ School of Statistics and Data Science, Shanghai University of }\\
{\small International Business and Economics, Shanghai, 201620, P.R. China} }

\date{}
\maketitle

\vspace{-0.5cm}

\begin{abstract}
A subset $D$ of the vertex set of a graph is called an $(a,b)$-regular set if each
vertex in $D$ has exactly $a$ neighbors in $D$, and each vertex outside $D$
has exactly $b$ neighbors in $D$.
A proper subgroup $H$ of $G$ is called an  $(a,b)$-regular set of $(G,\alpha)$ if it is an $(a,b)$-regular set of $GC(G,S,\alpha)$ for some generalized Cayley
subset $S$ of $G$ induced by the involutory automorphism $\alpha$.
This paper investigates subgroup $(0,b)$-regular sets in generalized Cayley graphs.
We first obtain some necessary and sufficient conditions for a proper subgroup of a group $G$ to be
a $(0,b)$-regular set of a generalized Cayley graph. In particular, we show that
a proper normal subgroup $H$ of $G$ is a $(0,b)$-regular set of $(G,\alpha)$ if and only if $H$ is a perfect code of $(G,\alpha)$ for some odd integer $b,$
which generalizes the result in [F.S. Seiedali, Z. Akhlaghi, B. Khosravi. Generalized Cayley graphs and perfect code. Discrete Math. 2026].
Furthermore,  for connected generalized Cayley graphs with $\alpha$-normal generalized
Cayley subsets, we  characterize subgroup $(0,b)$-regular sets
via equitable partitions.
\end{abstract}



\section{Introduction}
All groups considered in this paper are finite, and all graphs are finite and simple.
Let $\Gamma$ be a graph with vertex set $V(\Gamma)$ and edge set $E(\Gamma)$. For a vertex $v\in V(\Gamma)$,  denote its neighbourhood by $N_\Gamma(v)$. For a subset $X\subseteq V(\Gamma)$, the induced subgraph on $X$ is written as $\Gamma[X]$.
Let $a$ and $b$ be nonnegative integers.
A nonempty proper subset $D\subseteq V(\Gamma)$ is called an {\it $(a,b)$-regular set} if every vertex in $D$ has exactly $a$ neighbours in $D$ and  every vertex in $V(\Gamma)\setminus D$ has exactly $b$ neighbours in $D$. This notion unifies several classical concepts. $(0,1)$-regular sets correspond to perfect codes, and $(1,1)$-regular sets correspond to total perfect codes.
Perfect codes and total perfect codes have been extensively studied in Cayley graphs \cite{ChenZhou2016,ChenLiZhang2025}, Cayley sum graphs \cite{MaFengWang2020,Zhang2024,AhanjidehAkhlaghiRodrigues2026} and proper intersection power graphs \cite{M.L.Z}.
More generally, regular sets in these graph classes have attracted considerable attention in recent years, see \cite{KhaefiAkhlaghiKhosravi2025,SeiedaliKhosraviAkhlaghi2025,WangXu2024,WangXuZhou2024,ZhangYuan2026}.

Generalized Cayley graphs, as a generalization of Cayley graphs, were introduced by Maru\v{s}i\v{c}, Scapellato and Zagaglia Salvi \cite{Marusic1992}.
\begin{definition}\label{def1}
Let $G$ be a group,  $S\subseteq G$ and $\alpha \in \mathrm{Aut}(G)$ satisfy the following conditions
\begin{enumerate}[{\rm(a)}]

\item $\alpha^2=\mathrm{id}$, where $\mathrm{id}$ is the identity of $\mathrm{Aut}(G)$,

\item for any $g \in G$, $\alpha(g^{-1}) g\notin S$,

\item for $g,h\in G$, if $\alpha(h^{-1}) g\in S$, then $\alpha(g^{-1}) h\in S$.
\end{enumerate}
Then the graph with vertex set $G$ and edge set $\{\{g,h\}\mid \alpha(g^{-1}) h \in S\}$ is called a generalized Cayley graph and denoted by $GC(G,S,\alpha)$.
Here, $S$ is referred to as a generalized Cayley subset  of $G$ induced by $\alpha$, and $GC(G,S,\alpha)$ is the generalized Cayley graph with respect to the pair $(S, \alpha)$.
In particular, if $\alpha=\mathrm{id}$, then $GC(G, S, \alpha)$ is the Cayley graph $Cay(G,S)$.
\end{definition}

If a proper subgroup $H$ of a group $G$ is an $(a,b)$-regular set of $GC(G,S,\alpha)$,
then $H$ is called a {\it subgroup $(a,b)$-regular set} of $GC(G,S,\alpha)$.
A proper subgroup $H$ of $G$ is called an {\it $(a,b)$-regular set} of $(G,\alpha)$ if it is an $(a,b)$-regular set of $GC(G,S,\alpha)$ for some generalized Cayley
subset $S$ of $G$ induced by the involutory automorphism $\alpha$.
In recent years, generalized Cayley graphs have attracted increasing attention.
Their structural properties, including the isomorphism problem and connectivity,
were studied in \cite{H.K.P2017,ZhuLiuYang2023}.
In \cite{LiaoLiuLiTian2026}, Liao et al.\ established some necessary and sufficient
conditions for a subset of a group to be a perfect code or a total perfect code
in generalized Cayley graphs.
Seiedali et al.\ \cite{SeiedaliAkhlaghiKhosravi2026} also determined all
integer pairs $(a,b)$ for which a subgroup $H$ of an abelian group $G$ is an
$(a,b)$-regular set in generalized Cayley graphs of $G$.
These results indicate that classical notions in Cayley graphs and Cayley sum graphs,
such as perfect codes and regular sets, admit natural extensions to generalized
Cayley graphs.


The paper is organized as follows.
In Section 2, we present some notations and useful
lemmas.
In Section 3,  we obtain conditions under which a proper subgroup of a group is a $(0,b)$-regular set of a generalized Cayley graph.
In Section 4, for connected generalized Cayley graphs with $\alpha$-normal
generalized Cayley subsets,
we further characterize, via equitable partitions, when a proper subgroup is a $(0,b)$-regular set.

\section{Preliminaries}
Let $G$ be a group with identity element $e$, and let $\alpha$ be an involutory automorphism of $G$. Define  $\omega_\alpha(G)=\{\alpha(g^{-1})g\mid g\in G\}$, $\Omega_\alpha(G)=\{g\in G\setminus \omega_\alpha(G) \mid \alpha(g)=g^{-1}\}$ and $\mho_\alpha(G)=\{g\in G \mid \alpha(g)\neq g^{-1}\}$.
Set $T_\alpha(G)=\omega_\alpha(G)\cup \Omega_\alpha(G)$.
It is straightforward to check that $\alpha(x)=x^{-1}$ for every $x\in T_\alpha(G)$.
Moreover, the sets $\omega_\alpha(G)$, $\Omega_\alpha(G)$ and $\mho_\alpha(G)$ form a partition of $G$.
By Definition \ref{def1}, a subset $S$ of $G$ is a generalized Cayley subset induced by $\alpha$ if and only if it satisfies $S\cap \omega_\alpha(G)=\emptyset$ and $\alpha(S)=S^{-1}$.
A generalized Cayley subset $S$ of $G$ induced by $\alpha$  is called {\it $\alpha$-normal} if $\alpha(g^{-1})Sg=S$ holds for all $g\in G$.
For $x\in G$ and   a proper normal subgroup $H$  of $G$,
define $\mathcal{L}_x=|Hx-\omega_\alpha(G)|$
and  $\mathcal{L}(H)=\min\{\mathcal{L}_x\mid x\in G\setminus H\}$.

Throughout this section, $\alpha$ is an involutory automorphism of $G$, and $S$ is a generalized Cayley subset induced by $\alpha$. For $X\subseteq G$, write $\overline X=\sum_{x\in X}x$ in the integral group ring $\mathbb ZG$.

\begin{lemma}\label{lem1}
Let $G$ be a group and  $D$  a  nonempty proper subset of $G$. Let $a,b$ be nonnegative integers. Then the following statements are equivalent.
\begin{enumerate}[{\rm(a)}]
\item $D$ is an $(a,b)$-regular set of the generalized Cayley graph $GC(G,S,\alpha)$.

\item $|\alpha(x)S\cap D|=a$ for any $x\in D$, $|\alpha(x)S\cap D|=b$ for any $x\in G\setminus D$.

\item $\overline{\alpha(D)}\cdot \overline{S}=a \overline{D}+b \overline{G\setminus D}$.

\item  $\overline{\alpha(D)}\cdot \overline{S}+(b-a)\overline{D}=b \overline{G}$.
\end{enumerate}
\end{lemma}
\begin{proof}
It is straightforward to verify that conditions (a) and (b) are equivalent, and that conditions (c) and (d) are equivalent.
Next, we show that (b) and (c) are equivalent.
Since $\alpha(S^{-1})=S$, we have
\begin{align*}
\overline{\alpha(D)} \cdot \overline{S}
&= \left(\sum_{d \in D} \alpha(d)\right)\left(\sum_{s \in S} s\right) = \sum_{d \in D}\sum_{s \in S} \alpha(d)s \\
&= \sum_{x \in G}\left(\sum_{\substack{(d,s)\in D\times S \\ \alpha(d)s = x}} 1\right)x = \sum_{x \in G}\left(\sum_{\substack{d \in D \\ d \in \alpha(x)S}} 1\right)x \\
&= \sum_{x \in G} |D \cap \alpha(x)S|\, x = \sum_{x \in D} |D \cap \alpha(x)S|\, x + \sum_{x \in G\setminus D} |D \cap \alpha(x)S|\, x \\
&= a \, \overline{D} + b \, \overline{G \setminus D}.
\end{align*}

Hence, conditions (b) and (c) are equivalent.
\end{proof}

Based on Lemma \ref{lem1}, we obtain the following special cases for perfect codes and total perfect codes.

\begin{corollary}\label{cor1}
Let $G$ be a group and $D$  a  nonempty proper subset of $G$. Then
\begin{enumerate}[{\rm(a)}]
\item $D$ is a perfect code of the generalized Cayley graph $GC(G,S,\alpha)$ if and only if $\overline{\alpha(D)} \cdot \overline{S} + \overline{D} = \overline{G}.$

\item $D$ is a total perfect code of the generalized Cayley graph $GC(G,S,\alpha)$ if and only if $\overline{\alpha(D)} \cdot \overline{S} = \overline{G}.$
\end{enumerate}
\end{corollary}

Corollary \ref{cor1} provides  necessary and sufficient conditions for perfect codes and total perfect codes in generalized Cayley graphs, respectively, as established in \cite{LiaoLiuLiTian2026}.
Combining Lemma \ref{lem1} and Corollary \ref{cor1}, we obtain a relationship between perfect codes and $(0,b)$-regular sets in generalized Cayley graphs.

\begin{corollary}\label{cor2}
Let $G$ be a group and  $H$  a subgroup of $G$.
Let $b$ be a positive integer.
Assume that $H$ is a perfect code of the generalized Cayley graph $GC(G,S_1,\alpha)$.
Then $H$ is a $(0,b)$-regular set of the generalized Cayley graph  $GC(G,S_2,\alpha)$ if and only if $\overline{\alpha(H)} \cdot \overline{S_2}
= b \big(\overline{\alpha(H)} \cdot \overline{S_1}\big)$.
\end{corollary}

\begin{lemma}\label{lem2}
Let $G$ be a group and $S$  an $\alpha$-normal generalized Cayley subset of $G$. Let $b$ be a positive integer.
Then a subset $D \subseteq G$ is a $(0,b)$-regular set of $GC(G,S,\alpha)$ if and only if $Dg$ is a $(0,b)$-regular set of $GC(G,S,\alpha)$ for every $g \in G$.
\end{lemma}
\begin{proof}
Let $g \in G$. Since $S$ is an $\alpha$-normal generalized Cayley subset of $G$, the right translation $R_g: x \mapsto xg$ is an automorphism of $GC(G,S,\alpha)$, which  maps $D$ to $Dg$.
It follows that $D$ is a $(0,b)$-regular set of $GC(G,S,\alpha)$ if and only if $Dg$ is a $(0,b)$-regular set of $GC(G,S,\alpha)$.
\end{proof}

In the sequel, using a proof similar to that of ~\cite[Lemma~2.2]{WangXuZhou2024}, we have the following result for generalized Cayley graphs.
\begin{lemma}\label{lem3}
Let $b$ be a positive integer.
If $D$ is a $(0,b)$-regular set of the generalized Cayley graph $GC(G,S,\alpha)$, then $|G| = |D|\left(1 + \frac{|S|}{b}\right)$.
Moreover, if $H$ is a subgroup  $(0,b)$-regular set of $GC(G,S,\alpha)$, then $b$ divides $|S|$.
\end{lemma}


\section{Subgroup $(0,b)$-regular set}

In this section,
firstly,  we characterize when a proper subgroup of a group is a
$(0,b)$-regular set of a generalized Cayley graph, in terms of coset decompositions.

\begin{lemma}\label{lem3.1}
Let $H$ be a subgroup of a group $G$ and $b$ a positive integer.
If $H$ is a $(0,b)$-regular set of $GC(G,S,\alpha)$, then $\alpha(H)=H$ and $H\cap S=\emptyset$.
Moreover, there exists a subset $L\subseteq S$ such that $L\cup \{e\}$ forms a  set of right coset representatives of $H$ in $G$.
\end{lemma}
\begin{proof}
Since $H$ is a $(0,b)$-regular set of $GC(G,S,\alpha)$,
by Lemma \ref{lem1}, $|H\cap S|=0$, which implies $H\cap S=\emptyset$.
And since every $g\in G\setminus H$ is adjacent to exactly $b$ vertices in $H$, we have $
G\setminus H \subseteq \bigcup_{s\in S} \alpha(H)s.$
Note that for any $s_1,s_2\in S$, either $\alpha(H)s_1=\alpha(H)s_2$ or $\alpha(H)s_1\cap \alpha(H)s_2=\emptyset$.
Thus,  there exists a subset $L=\{\ell_1,\ell_2,\dots,\ell_r\}\subseteq S$ such that
\begin{equation*}
G\setminus H\subseteq\alpha(H)\ell_1\cup \alpha(H)\ell_2\cup \cdots \cup \alpha(H)\ell_r,
\end{equation*}
where $\alpha(H)\ell_i\cap \alpha(H)\ell_j=\emptyset$ for any $1\leq i\neq j\leq r.$
This yields
\begin{equation*}
G=H\cup \alpha(H)\ell_1\cup \alpha(H)\ell_2\cup \cdots \cup \alpha(H)\ell_r.
\end{equation*}

If $\alpha(H)=\alpha(H)\ell_i$ for some $1\leq i\leq r$, then there exists $h\in H$ such that $e=\alpha(h)\ell_i$, which implies that $e$ is adjacent to $h$, a contradiction to the assumption that $H$ is a $(0,b)$-regular set.
Hence $\alpha(H)\cap\alpha(H)\ell_i=\emptyset$ for every $1\leq i\leq r$.
Since \[ G=H\cup\alpha(H)\ell_1\cup\cdots\cup\alpha(H)\ell_r, \] it follows that $\alpha(H)\subseteq H$.
As $|\alpha(H)|=|H|$, we obtain $\alpha(H)=H$.
Consequently, \[ G=H\cup H\ell_1\cup\cdots\cup H\ell_r \] and $L\cup\{e\}$ forms a set of right coset representatives of $H$ in $G$.
\end{proof}

\begin{lemma}\label{lem3.2}
Let $H$ be a subgroup of a group $G$ and $b$ a positive integer.
If $H$ is a $(0,b)$-regular set of $GC(G,S,\alpha)$, then for each $s\in S$, there exist exactly $b-1$ distinct elements $s_1,s_2,\dots,s_{b-1}\in S\setminus \{s\}$ such that
\begin{equation*}
Hs_1=Hs_2=\cdots=Hs_{b-1}=Hs.
\end{equation*}
\end{lemma}
\begin{proof}
On the one hand,
by Lemma \ref{lem3.1}, we have $\alpha(H)=H$ and $H\cap S=\emptyset$.
Fix $s\in S$. Since $H$ is a $(0,b)$-regular set and $e$ is adjacent
to $s$, the vertex $s$ has exactly $b-1$ neighbors
$h_1,\ldots,h_{b-1}$ in $H\setminus\{e\}$. For each
$1\le i\le b-1$, there is a unique $s_i\in S$ such that $s=\alpha(h_i)s_i$.
Clearly, $s_i\neq s$ for all $1\le i\le b-1$, and the elements
$s_1,\ldots,s_{b-1}$ are pairwise distinct.
Since $\alpha(H)=H$, we
have $Hs_i=Hs$ for each $1\le i\le b-1$.

On the other hand,
let $t\in S\setminus\{s\}$ satisfy $Ht=Hs$.
We have $s=ht$ for some $h\in H\setminus\{e\}$ and $h=\alpha(h')$ for some $h'\in H\setminus\{e\}$
as $\alpha(H)=H$.
Hence $s=\alpha(h')t$, which implies that $h'$ is a neighbor of $s$ in $H\setminus\{e\}$. We have $h'=h_i$ for some $i$ and thus $t=s_i$.
Therefore, there are exactly $b-1$ distinct elements
$s_1,\ldots,s_{b-1}\in S\setminus\{s\}$ such that $Hs_1=\cdots=Hs_{b-1}=Hs$.
\end{proof}

\begin{theorem}\label{thm3.3}
Let $H$ be a proper subgroup of a  group $G$ with  $|G:H|=n$ and  $S$  a generalized Cayley subset of $G$  induced by  an involutory automorphism $\alpha$.
Let $b$ be a positive integer.
Then the following statements are equivalent.
\begin{enumerate}[{\rm(a)}]
\item $H$ is a $(0,b)$-regular set of $GC(G,S,\alpha)$.

\item $\alpha(H)=H$ and there exists a subset $L=\{\ell_1,\ell_2,\ldots,\ell_{n-1}\}$ of $S$ such that $L\cup \{e\}$ is a set of right coset representatives of $H$ in $G$, and $S$ can be partitioned into subsets $S_1,S_2,\ldots, S_{n-1}$ such that $|S_i|=b$ and $S\cap H\ell_i=S_i$ for $1\leq i\leq n-1$.

\item $\alpha(H)=H$ and $S$ can be partitioned into $b$ subsets $T_1,T_2,\ldots, T_b$ such that $|T_j|=n-1$ for each $1\le j\le b$, and $T_j\cup \{e\}$ is a set of right coset representatives of $H$ in $G$.
\end{enumerate}
\end{theorem}
\begin{proof}
(a)$\Rightarrow$ (b).
Since $H$ is a $(0,b)$-regular set of $GC(G,S,\alpha)$,
by Lemma \ref{lem3.1}, we have $\alpha(H)=H$, $S\cap H=\emptyset$ and there exists a subset $L=\{\ell_1,\dots,\ell_{n-1}\}\subseteq S$ such that $L\cup\{e\}$ is a  set of right coset representatives of $H$ in $G$.
Hence  there exists $\ell_i\in L$ such that $Hs=H\ell_i$ for each $s\in S$.
By Lemma \ref{lem3.2}, each $\ell_i$ corresponds to exactly $b-1$ elements in $S\setminus \{\ell_i\}$ lying in the same right coset, say $\{\ell_{i,1},\dots,\ell_{i,b-1}\}$.
Define $S_i=\{\ell_i,\ell_{i,1},\dots,\ell_{i,b-1}\}$.
Then $|S_i|=b$ and $\{S_1,\dots,S_{n-1}\}$ is a partition of $S$.

(b)$\Rightarrow$ (c).
Assume that (b) holds. Write $S_i=\{\ell_i,\ell_{i,1},\dots,\ell_{i,b-1}\}$ for any $1\leq i\leq n-1$.
Define $T_j=\{\ell_{1,j},\ell_{2,j},\dots,\ell_{n-1,j}\}$  for each $1\le j\le b-1$ and $T_b=\{\ell_1,\dots,\ell_{n-1}\}$.
Then $\{T_1,\dots,T_b\}$ is a partition of $S$, and each $T_j\cup\{e\}$ is a set of right coset representatives of $H$ in $G$.

(c)$\Rightarrow$ (a).
Assume that (c) holds.
Let $u\in G\setminus H$. Since
$T_j\cup\{e\}$ is a set of right coset representatives of $H$ in
$G$ for each $1\le j\le b$, there exists a unique
$s_j(u)\in T_j$ such that $Hs_j(u)=Hu$.
Hence there exists a unique $h_j(u)\in H$ such that $u=h_j(u)s_j(u)$.
Since  $\alpha^2=\mathrm{id}$, we have
$u=\alpha\bigl(\alpha(h_j(u))\bigr)s_j(u)$, which implies that $\alpha(h_j(u))$ is adjacent to $u$.
In fact,
these $b$ neighbors are pairwise distinct.
Otherwise, if $\alpha(h_i(u))=\alpha(h_j(u))$ for some $i\ne j$, then
$h_i(u)=h_j(u)$, and thus
\[
s_i(u)=h_i(u)^{-1}u=h_j(u)^{-1}u=s_j(u),
\]
contradicting  the fact that $T_i$ and $T_j$ are disjoint.
Thus, since $\alpha(H)=H$,  $u$ has at least $b$ distinct neighbors in $H$.

It remains to show that $\{\alpha(h_j(u))\mid 1\leq j\leq b\}$ are all the neighbors of $u$ in $H$.
Let $h\in H$ be any neighbor of $u$.
Then $u=\alpha(h)s$ for some $s\in S$.
Since $S=\cup_{j=1}^b T_j$, there is a unique $j$ such that $s\in T_j$, where $1\leq j\leq b.$
As $\alpha(h)\in H$, we have $Hs=Hu$.
By the uniqueness of the representative of $Hu$ in $T_j$, it
follows that $s=s_j(u)$.
Hence $\alpha(h)=h_j(u)$, and thus $h=\alpha(h_j(u))$,
which implies that every neighbor of $u$ in $H$ is one of $\alpha(h_1(u)),\ldots,\alpha(h_b(u))$.

Now we prove that no two vertices in $H$ are adjacent. Otherwise,
there exists some vertex $h_{1}\in H$ such that $h_{1}$ is adjacent to $h_{2},$ where $h_{2}\in H.$ Then
there exists $s\in S$ such that  $\alpha(h_{1}^{-1})h_{2}=s.$ Since $S=\cup_{j=1}^b T_j$, there is a unique $j$ such that $s\in T_j$.
We have $Hs=H\alpha(h_{1}^{-1})h_{2}=H$ by $\alpha(H)=H,$
which contradicts  the fact that $T_j\cup \{e\}$ is a set of right coset representatives of $H$ in $G$.
Hence $H$ is a $(0,b)$-regular set of $GC(G,S,\alpha)$.
\end{proof}

To further interpret $(0,b)$-regular sets in a graph-theoretic framework, we recall the notion of graph coverings.
A graph $\Gamma$ is called a {\it $b$-cover} of a graph $\Sigma$ if there exists a surjective homomorphism $f\colon \Gamma\to \Sigma$ such that $|N_\Gamma(u)\cap f^{-1}(y)|=b$ for every $\{x,y\}\in E(\Sigma)$ and all $u\in f^{-1}(x)$.
And the corresponding homomorphism $f$ is called a {\it $b$-covering} from $\Gamma$ to $\Sigma$.
In particular, if $f$ is $m$-to-$1$ for some integer $m\ge b$, then $f$ is called an {\it $m$-fold $b$-covering} and the graph $\Gamma$ is called an {\it $m$-fold $b$-cover} of $\Sigma$.
Moreover, for $x\in V(\Sigma)$, the set $f^{-1}(x)$ is called the {\it fibre} of $x$, and
\[
f^{-1}(\{x,y\})=\{\{u,v\}\in E(\Gamma)\mid \{f(u),f(v)\}=\{x,y\}\}
\]
is called the {\it fibre} of the edge $\{x,y\}$.

We now summarize several basic properties of a $b$-covering $f\colon \Gamma\to \Sigma$.
Firstly, the fibre of each vertex in $\Sigma$ is an independent set.
Secondly, for any two vertices $u,v\in V(\Sigma)$, if $u$ is adjacent to $v$, then the induced subgraph $\Gamma[f^{-1}(u)\cup f^{-1}(v)]$ is a $b$-regular bipartite graph. If $u$ is  not adjacent to $v$, then the induced subgraph $\Gamma[f^{-1}(u)\cup f^{-1}(v)]$ is an empty graph.
Furthermore, for any $\{x,y\}\in E(\Sigma)$, counting edges between the vertex-fibres $f^{-1}(x)$ and $f^{-1}(y)$ shows that $|f^{-1}(x)|=|f^{-1}(y)|$.
Hence, if $\Sigma$ is connected and has at least two vertices, then all vertex-fibres have equal size, and $f$ is an $m$-fold $b$-covering for some integer $m$.

\begin{example}
The graphs $C_6$ and $K_2$ are shown in Figure \ref{fig1}. Define $f:u_1,u_2,u_3\mapsto w_1$ and $v_1,v_2,v_3\mapsto w_2$. Then $f$ is a $3$-fold $2$-covering from $C_6$ to $K_2$.
\begin{figure}[h]
\centering
\includegraphics[scale=0.6]{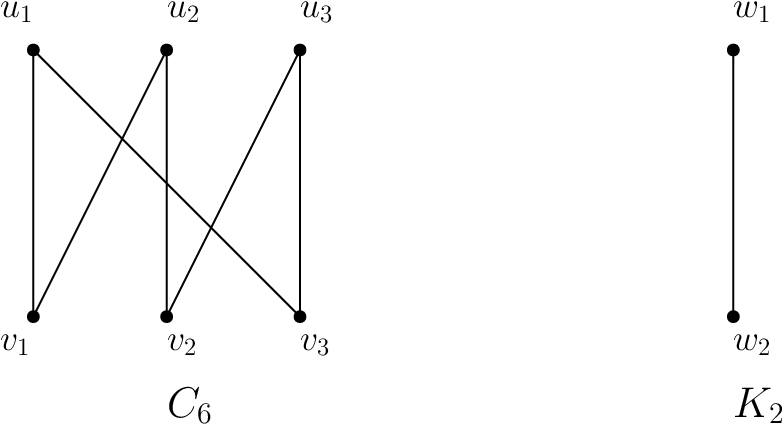}
\caption{$C_6$ and  $K_2$}
\label{fig1}
\end{figure}
\end{example}

The following lemma characterizes the relationship between $(0,b)$-regular sets and $b$-coverings.
Since the proof was omitted in \cite{WangXuZhou2024}, we include a concise proof here.

\begin{lemma}[\cite{WangXuZhou2024}]\label{lem3.3}
Let $\Gamma$ and $\Sigma$ be two graphs, and let $b$ be a positive integer. Then the following hold.
\begin{enumerate}[{\rm(a)}]
\item If $f$ is a $b$-covering from $\Gamma$ to $\Sigma$, and $D$ is a $(0,1)$-regular set of $\Sigma$, then $f^{-1}(D)$ is a $(0,b)$-regular set of $\Gamma$.

\item If $\Gamma$ is a $d$-regular graph, where $d\geq 1,$ and $D_1,D_2,\dots,D_n$ are pairwise disjoint $(0,b)$-regular sets of $\Gamma$, then the induced subgraph $\Gamma[\bigcup_{\ell=1}^n D_\ell]$ is an $m$-fold $b$-cover of the complete graph $K_n$, where $m=|D_\ell|=\frac{b|V(\Gamma)|}{b+d}.$
\end{enumerate}
\end{lemma}
\begin{proof}
(a)  Assume that $f$ is a $b$-covering  from $\Gamma$ to $\Sigma$, and let $D$ be a $(0,1)$-regular set of $\Sigma$.
For any vertex $u\in V(\Gamma)\setminus f^{-1}(D)$, let $v=f(u)\notin D$.
Since $D$ is a $(0,1)$-regular set of $\Sigma$, the vertex $v$ has exactly one neighbor $x$ in $D$.
By the definition of  $b$-coverings, we have $|N_\Gamma(u)\cap f^{-1}(x)| = b.$
Since $f^{-1}(x)\subseteq f^{-1}(D)$, it follows that $u$ has exactly $b$ neighbors in $f^{-1}(D)$.
Let $w\in f^{-1}(D)$ and  $z=f(w)\in D$.
Since $z$ has no neighbors in $D$, it follows that $w$ has no neighbors in $f^{-1}(D)$.
Thus, $f^{-1}(D)$ is a $(0,b)$-regular set of $\Gamma$.

(b)
Let $\Gamma$ be a $d$-regular graph, and let $D_1,D_2,\dots,D_n$ be pairwise disjoint $(0,b)$-regular sets of $\Gamma$.
For each $D_i$, counting edges between $D_i$ and $V(\Gamma)\setminus D_i$ gives $d|D_i| = b\big(|V(\Gamma)|-|D_i|\big),$
which yields $|D_i|=\frac{b|V(\Gamma)|}{b+d}.$

We now show that the induced subgraph $\Gamma\left[\bigcup_{\ell=1}^n D_\ell\right]$ is an $m$-fold $b$-cover of $K_n$, where $m=|D_\ell|$.
Since the sets $D_1,\dots,D_n$ are pairwise disjoint $(0,b)$-regular sets, each induced subgraph $\Gamma[D_i\cup D_j]$ is a $b$-regular bipartite graph for any $i\neq j$.
Denote by  $V(K_{n})=\{v_1,\dots,v_n\}$.
Define a surjective mapping
\[
f\colon \Gamma\left[\bigcup_{\ell=1}^n D_\ell\right]\to K_n
\]
by sending every vertex in $D_i$ to $v_i$.
It is straightforward to verify that $f$ is a graph homomorphism.
Moreover, for any two adjacent vertices $v_i,v_j$ in $K_n$ and any vertex $u\in D_i$, we have $|N_{\Gamma[\cup_{\ell=1}^n D_\ell]}(u)\cap D_j|=b.$
Therefore, $\Gamma\left[\bigcup_{\ell=1}^n D_\ell\right]$ is an $m$-fold $b$-cover of $K_n$, where $m=|D_\ell|=\frac{b|V(\Gamma)|}{b+d}.$
\end{proof}

In particular,
for $\alpha$-normal generalized Cayley subsets, we   give a characterization of subgroup
$(0,b)$-regular sets
in terms of $b$-coverings.


\begin{theorem}\label{thm3.6}
Let $H$ be a proper subgroup of a group $G$ and $S$  an $\alpha$-normal generalized Cayley subset of $G.$ Let $b$ be a positive integer.
Then $H$ is a $(0,b)$-regular set of $GC(G,S,\alpha)$ if and only if there exists a $b$-covering  $f$ from $GC(G,S,\alpha)$ to $K_n$ such that $H$ is a vertex-fibre of $f$, where $n=|G:H|$.
\end{theorem}
\begin{proof}
Let $n=|G:H|$.
Suppose that $H$ is a $(0,b)$-regular set of $GC(G,S,\alpha)$.
By Lemma \ref{lem2}, every right coset $Hs$ with $s\in S$ is also a $(0,b)$-regular set of $GC(G,S,\alpha)$.
Moreover, by Theorem \ref{thm3.3}, there exists a subset $L=\{\ell_1,\ell_2, \ldots, \ell_{n-1}\}\subseteq S$ such that $L\cup\{e\}$ is a  set of right coset representatives of $H$ in $G$.
Thus,  $H, H\ell_1,H\ell_2, \ldots, H\ell_{n-1}$ are pairwise disjoint $(0,b)$-regular sets of $GC(G,S,\alpha).$
Applying Lemma \ref{lem3.3}(b), we obtain an $|H|$-fold $b$-covering
$f$ from  $GC(G,S,\alpha)$ to $K_n$
whose vertex-fibres are  $H$ and $H\ell_i$, where $1\leq i\leq n-1$.

Conversely, suppose that there exists a $b$-covering
$f$ from $ GC(G,S,\alpha)$ to $ K_n$
such that $H$ is a vertex fibre of $f$.
Write $H=f^{-1}(x)$ for some $x\in V(K_n)$.
Since $\{x\}$ is a $(0,1)$-regular set of $K_n$, Lemma \ref{lem3.3}(a) implies that $H$ is a $(0,b)$-regular set of $GC(G,S,\alpha)$.
\end{proof}

Next, we consider the conditions under which a proper normal subgroup of a group is a
$(0,b)$-regular set of a generalized Cayley graph.

\begin{theorem}\label{thm3.7}
Let $\alpha$ be an involutory automorphism of a group $G$ and $H$  a proper normal subgroup of $G$ such that $\alpha(H)=H$.
Let $b$ be a positive even integer with $b\leq \mathcal{L}(H)$.
Then $H$ is a $(0,b)$-regular set of $(G, \alpha)$.
\end{theorem}
\begin{proof}
Let $|G:H|=n$ and $L=\{e, a_1, \ldots, a_{n-1}\}$ be a set of right coset representatives of $H$ in $G$.
Assume $|Ha_i\cap \Omega_\alpha(G)|=m_i$ for $1\leq i\leq n-1$.
Let $\sigma$ be a permutation on $\{1,2,\ldots,n-1\}$ such that $H\alpha(a_i^{-1})=Ha_{\sigma(i)}$.
Since $H$ is a normal subgroup of $G$ and $\alpha(H)=H$, we have
\begin{equation*}
Ha_{\sigma^2(i)}=H\alpha(a_{\sigma(i)}^{-1})=\alpha(Ha_{\sigma(i)}^{-1})
=\alpha(Ha_{\sigma(i)})^{-1}=\alpha(H\alpha(a_i^{-1}))^{-1}=Ha_i.
\end{equation*}
It follows that $\sigma$ is the identity permutation or an involution.
Assume $\sigma$ fixes $t$ integers in $\{1,2,\ldots,n-1\}$.
Then $0\leq t\leq n-1$ and $n-1-t$ is even.
Let $\ell=\frac{n-1-t}{2}$.
Without loss of generality, we assume
\begin{equation*}
\sigma(i)=
\begin{cases}
i, \quad \quad \quad i\leq t,\\
i+\ell, \quad t<i\leq t+\ell,\\
i-\ell, \quad t+\ell<i\leq n-1.
\end{cases}
\end{equation*}
Then
\begin{equation*}
\alpha(Ha_i)^{-1}=H\alpha(a_i^{-1})=Ha_{\sigma(i)}=Ha_i ~\text{for}~ 0<i\leq t
\end{equation*}
and
\begin{equation*}
\alpha(Ha_{t+j})^{-1}=H\alpha(a_{t+j}^{-1})=Ha_{\sigma(t+j)}=Ha_{t+j+\ell} ~\text{for}~ 0<j\leq \ell.
\end{equation*}

If $i>t$, then $\sigma(i)\neq i$.
In this case, we first show that $Ha_i\cap T_\alpha(G)=\emptyset$.
Suppose, to the contrary, that there exists some
$x\in Ha_i\cap T_\alpha(G)$.
Since $x\in T_\alpha(G)$,
we have $\alpha(x^{-1})=x$.
Moreover, since $x\in Ha_i$ and
$\alpha(H)=H$, we have
\[
\alpha(x^{-1})
\in \alpha((Ha_i)^{-1})
=\alpha(a_i^{-1}H)
=\alpha(a_i^{-1})H
=H\alpha(a_i^{-1})
=Ha_{\sigma(i)}.
\]
Hence $x\in Ha_i\cap Ha_{\sigma(i)}$, which implies $Ha_i=Ha_{\sigma(i)}$,
a contradiction to $\sigma(i)\neq i$.
Therefore, $Ha_i\cap T_\alpha(G)=\emptyset$, which implies that
$Ha_i\subseteq \mho_\alpha(G)$  for any $i>t.$

For any even integer $b$ with $b\leq \mathcal{L}(H)$, now we  construct a subset $S_i\subseteq Ha_i$ of size $b$ for each $i$.
\begin{enumerate}[{\rm$\bullet$}]
\item When $i\le t$ and $m_i\ge b$, choose any $b$ elements of $Ha_i\cap\Omega_\alpha(G)$.

\item When $i\le t$, $m_i<b$ and $b-m_i$ is even, let $S_i$ be the disjoint union of $m_i$ elements of $Ha_i\cap\Omega_\alpha(G)$ and $\frac{b-m_i}{2}$ pairs $\{x,\alpha(x^{-1})\}$ from $Ha_i\cap\mho_\alpha(G)$.

\item When $i\le t$, $m_i<b$ and $b-m_i$ is odd, we have $m_i\ge1$.
      Let $S_i$ be the disjoint union of $m_i-1$ elements of $Ha_i\cap\Omega_\alpha(G)$ and $\frac{b+1-m_i}{2}$ pairs $\{x,\alpha(x^{-1})\}$ from $Ha_i\cap\mho_\alpha(G)$.

\item When $t<i\leq t+\ell$, choose any subset $S_i\subseteq Ha_i$ of size $b$.

\item When $i>t+\ell$,  define $S_i=\alpha(S_{i-\ell})^{-1}$.
\end{enumerate}
Let $S=\cup_{i=1}^{n-1} S_i$.
By the construction of $S_{i}$, we have $S\cap \omega_\alpha(G)=\emptyset$ and $\alpha(S^{-1})=S$,
which implies that  $S$ is a generalized Cayley subset of $G$ with respect to $\alpha$.
For each $1\leq i\leq n-1,$ choose $s_{i}\in S_{i}$.   Then $Hs_{i}=Ha_{i}.$ Thus,
$\{e,s_{1},\ldots,s_{n-1}\}$ is  a set of right coset representatives of $H$ in $G$.
 Since $S\cap Hs_{i}=S\cap Ha_{i}=S_{i}$ and $|S_{i}|=b$ for every $1\leq i\leq n-1,$  Theorem \ref{thm3.3} implies that  $H$ is a $(0,b)$-regular set of $GC(G,S,\alpha)$.
\end{proof}

\begin{lemma}[\cite{LiaoLiuLiTian2026}]\label{lem5.3}
Let $G$ be a group and $\alpha$ an involutory automorphism of $G$.
Let $H$ be a proper subgroup of $G$.
Then $H$ is a perfect code of $G$ with respect to $\alpha$ if and only if $\alpha(H)=H$, and for any $x\in G\setminus H$ with $x\alpha(x)\in H$ and $\frac{|H|}{|H^x\cap H|}$ being odd, there exists $y\in Hx$ such that $y\in \Omega_\alpha(G)$.
\end{lemma}

By Lemma \ref{lem5.3}, if $H$ is a normal subgroup of $G$, then $H$ is a  perfect code of $G$ with respect to $\alpha$ if and only if $\alpha(H)=H$, and for any $x\in G\setminus H$ with $x\alpha(x)\in H$, there exists $y\in Hx$ such that $y\in \Omega_\alpha(G)$.

\begin{theorem}\label{thm3.9}
Let $\alpha$ be an involutory automorphism of a group $G$ and  $H$  a proper normal subgroup of $G$.
Let $b$ be an odd integer with $1\leq b\leq|H|$.
If $H$ is a $(0,b)$-regular set of $(G, \alpha)$, then $H$ is a perfect code of $(G, \alpha)$.
\end{theorem}
\begin{proof}
Assume that $H$ is a $(0,b)$-regular set of $GC(G,S,\alpha)$.
By Theorem~\ref{thm3.3} (c), $\alpha(H)=H$ and $S$ admits a decomposition
$S=\bigcup_{i=1}^b T_i$,
where each $T_i\cup \{e\}$ is a set of right coset representatives of $H$ in $G$.
To apply Lemma~\ref{lem5.3}, it suffices to show that, for every $g\in G\setminus H$, if $g\alpha(g)\in H,$ then $Hg\cap \Omega_\alpha(G)\neq\emptyset$. Suppose  $g\in G\setminus H$ satisfies $g\alpha(g)\in H$.
For each $1\le i\le b$, there exists a unique element $s_{i,j_i}\in T_i$ such that $g\in Hs_{i,j_i}$.
Set $T=\{s_{i,j_i}\in S\mid Hg=Hs_{i,j_i}\}$.
Then  $T\subseteq S$ and  $|T|=b$.
Since $g\alpha(g)\in H$,
\begin{equation*}
Hs_{i,j_i}=Hg=H\alpha(g^{-1})=\alpha(Hg)^{-1}=\alpha(Hs_{i,j_i})^{-1}=H\alpha(s_{i,j_i})^{-1}.
\end{equation*}
Thus   $\alpha(s_{i,j_i})^{-1}\in T$ if $s_{i,j_i}\in T$, which implies that  $\alpha(T^{-1})=T$.
Observe that $S\cap \omega_\alpha(G)=\emptyset$, hence $T\subseteq \Omega_\alpha(G)\cup \mho_\alpha(G)$.
Since  $|T|=b$ is odd and $T\cap \mho_\alpha(G)$ is a union of pairs $\{x, \alpha(x^{-1})\}$, it follows that $T\cap \Omega_\alpha(G)\neq\emptyset$.
Thus, $Hg\cap \Omega_\alpha(G)\neq\emptyset$.
By Lemma \ref{lem5.3}, $H$ is a perfect code of $(G,\alpha)$.
\end{proof}

When $b$ is odd, Theorem~\ref{thm3.9} shows that being a perfect code of $(G,\alpha)$ is a necessary condition for a proper normal subgroup
$H$ of $G$
to be  a $(0,b)$-regular set of $(G,\alpha)$.
We now turn to  sufficient conditions under which $H$ becomes a $(0,b)$-regular set of $(G,\alpha)$ for an odd integer $b$.
Since the proof of the following theorem is essentially the same as that of
Theorem~\ref{thm3.7}, we   present only a shortened version here.

\begin{theorem}\label{thm3.10}
Let $\alpha$ be an involutory automorphism of a group $G$ and  $H$  a proper normal subgroup of $G$.
If $H$ is a perfect code  of $(G, \alpha)$, then $H$ is a $(0,b)$-regular set of $(G, \alpha)$ for each odd integer $b$ with $1\leq b\leq \mathcal{L}(H)$.
\end{theorem}
\begin{proof}
Assume $H$ is a perfect code  of $(G, \alpha)$.
By Lemma \ref{lem5.3}, $\alpha(H)=H$, and for any $x\in G\setminus H$ such that $x\alpha(x)\in H$, there exists $y\in Hx$ such that $y\in \Omega_\alpha(G)$.

Let $|G:H|=n$ and $L=\{e, a_1, \ldots, a_{n-1}\}$ be a set of right coset representatives of $H$ in $G$.
Assume $|Ha_i\cap \Omega_\alpha(G)|=m_i$ for $1\leq i\leq n-1$.
Let $\sigma$ be a permutation on $\{1,2,\ldots,n-1\}$ such that $H\alpha(a_i^{-1})=Ha_{\sigma(i)}$.
As in the proof of Theorem~\ref{thm3.7}, $\sigma$ is either the identity permutation or an involution and we   assume
\begin{equation*}
\sigma(i)=
\begin{cases}
i, \quad \quad \quad i\leq t,\\
i+\ell, \quad t<i\leq t+\ell,\\
i-\ell, \quad t+\ell<i\leq n-1,
\end{cases}
\end{equation*}
where $\ell=\frac{n-1-t}{2}$.
Then $H\alpha(a_i^{-1})=Ha_i ~\text{for}~ 0<i\leq t$
and
$H\alpha(a_{t+j}^{-1})=Ha_{t+j+\ell} ~\text{for}~ 0<j\leq \ell$.
For any $i\leq t$, the condition $H\alpha(a_i^{-1})=Ha_i$ implies $a_i\alpha(a_i)\in H$.
Since $H$ is a perfect code of $(G, \alpha)$, Lemma~\ref{lem5.3} ensures that
$Ha_i\cap\Omega_\alpha(G)\neq\emptyset$, and hence $m_i\ge1$.
For any $i>t$, the same argument as in the proof of Theorem~\ref{thm3.7} shows that $Ha_i\subseteq \mho_\alpha(G)$.

Let $b$ be an odd integer with $1\leq b\le\mathcal{L}(H)$.
For each coset $Ha_i$, we construct a subset $S_i\subseteq Ha_i$ of size $b$ as follows.
\begin{enumerate}[{\rm$\bullet$}]
\item When $i\le t$ and $m_i\ge b$, choose any $b$ elements of $Ha_i\cap\Omega_\alpha(G)$.

\item When $i\le t$, $1\le m_i<b$ and $b-m_i$ is even, let $S_i$ be the disjoint union of $m_i$ elements of $Ha_i\cap\Omega_\alpha(G)$ and $\frac{b-m_i}{2}$ pairs $\{x,\alpha(x^{-1})\}$ from $Ha_i\cap\mho_\alpha(G)$.

\item When $i\le t$, $1\le m_i<b$ and $b-m_i$ is odd, let $S_i$ be the disjoint union of $m_i-1$ elements of $Ha_i\cap\Omega_\alpha(G)$ and $\frac{b+1-m_i}{2}$ pairs $\{x,\alpha(x^{-1})\}$ from $Ha_i\cap\mho_\alpha(G)$.

\item When $t<i\leq t+\ell$, choose any subset $S_i\subseteq Ha_i$ of size $b$.

\item When $i>t+\ell$, define $S_i=\alpha(S_{i-\ell})^{-1}$.
\end{enumerate}
Let $S=\bigcup_{i=1}^{n-1} S_i$.
Then $S$ is a generalized Cayley subset of $G$ with respect to $\alpha$ satisfying the condition~(c) of Theorem~\ref{thm3.3}.
By Theorem~\ref{thm3.3}, it follows that $H$ is a $(0,b)$-regular set of $(G,\alpha)$.
\end{proof}

Based on  Theorems \ref{thm3.9} and \ref{thm3.10}, we obtain the following corollary.

\begin{corollary}\label{cor3.11}
Let $\alpha$ be an involutory automorphism of a group $G$ and  $H$  a proper normal subgroup of $G$.
Let $b$ be an odd integer with $1\leq b\leq\mathcal{L}(H)$.
Then $H$ is a $(0,b)$-regular set of $(G, \alpha)$ if and only if $H$ is a perfect code of $(G, \alpha)$.
\end{corollary}

\begin{remark}
When $G$ is abelian, Theorem~\ref{thm3.7} and Corollary \ref{cor3.11} reduce to Lemma~2.11 in~\cite{SeiedaliAkhlaghiKhosravi2026}.
Thus, Lemma~2.11 in~\cite{SeiedaliAkhlaghiKhosravi2026} can be viewed as a special case of our results.
\end{remark}

Finally, we   establish  a relation between subgroup $(0,b)$-regular sets and
subgroup $(a,b)$-regular sets for certain values of $a$.
Let $H$ be a subgroup of a group $G,$ and let $\alpha$ be an involutory automorphism of $G$. Define
\[
\mathcal{A}_{\alpha}(H)
=
\left\{
|T| \mid
T\subseteq H-\omega_{\alpha}(G)
\text{ and }
\alpha(T^{-1})=T
\right\}.
\]

\begin{theorem}\label{thm3.13}
Let $\alpha$ be an involutory automorphism of a group $G$ and
$H$  a proper subgroup of $G$.
Let $a$ be a nonnegative integer with
$a\in\mathcal{A}_{\alpha}(H)$ and $b$ a positive integer.
Then $H$ is a $(0,b)$-regular set of
$(G,\alpha)$ if and only if $\alpha(H)=H$ and $H$ is an
$(a,b)$-regular set of $(G,\alpha)$.
\end{theorem}

\begin{proof}
Firstly, suppose  that $\alpha(H)=H$ and  $H$ is an
$(a,b)$-regular set of $GC(G,S,\alpha)$ for some generalized Cayley
subset $S$ of $G$ induced by $\alpha$.
By Lemma~2.1,
$|H\cap \alpha(g)S|=a$ for every $g\in  H$ and $|H\cap\alpha(g)S|=b$ for every $g\in G\setminus H$.
Set $T=S-H$.
Since $S\cap\omega_{\alpha}(G)=\emptyset$, we have
$T\cap\omega_{\alpha}(G)=\emptyset$.
Moreover, since $\alpha(H)=H$,
\begin{equation*}
\alpha(T^{-1})=\alpha\bigl((S-H)^{-1}\bigr)
=\alpha(S^{-1})-\alpha(H^{-1})
=S-H=T.
\end{equation*}
Thus $T$ is a generalized Cayley subset of $G$ induced by $\alpha$.

Observe that $|H\cap \alpha(g)T|=|\alpha(g^{-1})H\cap T|=|H\cap T|=0$ for every $g\in  H$.
Since $\alpha(g)\notin H$ for every $g\in G\setminus H$, the coset
$\alpha(g)H$ is disjoint from $H$. Hence
\[
H\cap\alpha(g)(S\cap H)=\emptyset,
\]
and consequently
\[
|H\cap\alpha(g)T|
=
|H\cap\alpha(g)S|
-|H\cap\alpha(g)(S\cap H)|
=|H\cap\alpha(g)S|
=
b.
\]
By Lemma~2.1, $H$ is a $(0,b)$-regular set of
$GC(G,T,\alpha)$. Therefore, $H$ is a $(0,b)$-regular set of
$(G,\alpha)$.

Conversely, suppose that $H$ is a $(0,b)$-regular set of
$GC(G,T,\alpha)$ for some generalized Cayley subset $T$ of $G$
induced by $\alpha$. By Lemma~3.1, $\alpha(H)=H$ and $H\cap T=\emptyset$.
Moreover, by Lemma~2.1, $|H\cap\alpha(g)T|=0$ for every $g\in  H$ and $|H\cap\alpha(g)T|=b$ for every $g\in G\setminus H$.

Since $a\in\mathcal{A}_{\alpha}(H)$, there exists a subset
$T'\subseteq H-\omega_{\alpha}(G)$ such that
$|T'|=a$ and $\alpha((T')^{-1})=T'$.
Set $S=T\cup T'$. Then $S\cap\omega_{\alpha}(G)=\emptyset$.
Since $T\cap H=\emptyset$ and $T'\subseteq H$, the union $T\cup T'$ is  disjoint, which implies that
\[
\alpha(S^{-1})
=
\alpha(T^{-1})\cup\alpha((T')^{-1})
=
T\cup T'
=
S.
\]
Hence $S$ is a generalized Cayley subset of $G$ induced by $\alpha$.
Note that  $|H\cap S|=|H\cap T'|=|T'|=a$.
Thus, for every $g\in H$,
\begin{equation*}
|H\cap \alpha(g)S|=|\alpha(g^{-1})H\cap S|=|H\cap S|=a.
\end{equation*}
Since $\alpha(H)= H,$  we have $\alpha(g)\notin H$ for every $g\in G\setminus H$, which implies that the coset $\alpha(g)H$ is disjoint from $H$.
Consequently,  $H\cap\alpha(g)T'\subseteq H\cap\alpha(g)H =\emptyset$.
Therefore,
\[
|H\cap\alpha(g)S|
=
|H\cap\alpha(g)T|
=
b.
\]
By Lemma~2.1, $H$ is an $(a,b)$-regular set of
$GC(G,S,\alpha)$. Thus $H$ is an $(a,b)$-regular set of
$(G,\alpha)$.
\end{proof}

\section{$(0,b)$-regular sets in connected generalized Cayley graphs}

In this section, we investigate $(0,b)$-regular sets in connected generalized Cayley graphs using equitable partitions.

Let $\Gamma$ be a graph with vertex set $V(\Gamma)$ and edge set $E(\Gamma)$, and let $A(\Gamma)$ denote its adjacency matrix.
A partition $\pi=\{V_1,V_2,\dots,V_r\}$ of $V(\Gamma)$ is called an {\it equitable partition} if there exists an $r\times r$ matrix $M=(m_{ij})$, called the {\it quotient matrix } of $\pi$, such that every vertex in $V_i$ has exactly $m_{ij}$ neighbors in $V_j$ for all $1\le i,j\le r$.

Let $\Gamma$ be a connected $d$-regular graph.
An equitable partition $\pi$ is called a {\it $\mu$-equitable partition} if all eigenvalues of its quotient matrix $M$, except the principal eigenvalue $d$, are equal to $\mu$.
A nonempty proper subset $D\subseteq V(\Gamma)$ is called a {\it $\mu$-perfect set} of $\Gamma$ if $\{D, V(\Gamma)\setminus D\}$ forms a $\mu$-equitable partition of $\Gamma$.

\begin{lemma}[\cite{BaileyCameron2019}]\label{lem4.1}
Let $\Gamma$ be a connected regular graph and  $\pi=\{V_1,V_2,\ldots, V_r\}$  a partition of $V(\Gamma)$ with  $r\geq 2$.
If $\pi$ is $\mu$-equitable, then each $V_i$ is a $\mu$-perfect set of $\Gamma$.
Conversely, if every $V_i$ is a $\mu$-perfect set of $\Gamma$, then $\pi$ is $\mu$-equitable.
\end{lemma}

\begin{lemma}[\cite{WangXuZhou2024}]\label{lem4.2}
Let $\Gamma$ be a connected regular graph and  $\pi=\{V_1,V_2,\dots,V_r\}$  a $(-b)$-equitable partition of $V(\Gamma)$ with  $r\geq 2$.
If for some $i$, $V_i$ is an independent set of $\Gamma$ , then $b$ is a positive integer and $V_i$ is a $(0,b)$-regular set of  $\Gamma$.
\end{lemma}

\begin{theorem}\label{thm4.3}
Let $H$ be a proper subgroup of a group $G$ and $S$  an $\alpha$-normal generalized Cayley subset of $G$ such that $GC(G,S,\alpha)$ is connected.
Let $b$ be an integer with $1\leq b\leq |S|$.
Then the following statements are equivalent.
\begin{enumerate}[{\rm(a)}]
\item $H$ is a $(0,b)$-regular set of $GC(G,S,\alpha)$.

\item $\alpha(H)=H$, $H\cap S=\emptyset$ and $GC(G,S,\alpha)$ admits a $(-b)$-equitable partition into $|G:H|$ blocks, one of which is $H$.

\item For each integer $0\le a\le \frac{|S|}{b}-1$, there exists an $(ab,(a+1)b)$-regular set of $GC(G,S,\alpha)$ consisting of $a+1$ right cosets of $H$, one of which is $H$.
\end{enumerate}
\end{theorem}
\begin{proof}
Let $|G:H|=n$.
We establish the equivalence by proving $(a)\Leftrightarrow (b)$ and $(a)\Leftrightarrow (c)$.

(a)$\Rightarrow$(b).
Since $H$ is a $(0,b)$-regular set of $GC(G,S,\alpha)$,
by Lemma \ref{lem3.1}, we have $\alpha(H)=H$ and $H\cap S=\emptyset$. Moreover, there exists a subset $L=\{\ell_1,\ell_2,\dots,\ell_{n-1}\}\subseteq S$ such that $L\cup \{e\}$ is a set of right coset representatives of $H$ in $G$.
By Lemma \ref{lem2}, each $H\ell_i$ is also a $(0,b)$-regular set of $GC(G,S,\alpha)$. Thus,
for each right coset $H\ell_i$, the induced partition $\{H\ell_i, G\setminus H\ell_i\}$ has quotient matrix
\[
\begin{pmatrix}
0 & |S|\\
b & |S|-b
\end{pmatrix},
\]
whose non-principal eigenvalue is $-b$.
Hence each $H\ell_i$ is a $(-b)$-perfect set of $GC(G,S,\alpha)$.
Applying Lemma \ref{lem4.1}, we obtain that $\{H, H\ell_1,\dots,H\ell_{n-1}\}$ forms a $(-b)$-equitable partition of $GC(G,S,\alpha)$.

(b)$\Rightarrow$(a).
If there exist $x,y\in H$ such that $x$ is adjacent to $y$, then $\alpha(x^{-1})y\in S$.
Since $\alpha(H)=H$, we obtain $\alpha(x^{-1})y\in H\cap S$, which contradicts the assumption that $H\cap S=\emptyset$.
Thus $H$ is an  independent set of $GC(G,S,\alpha)$.
By Lemma \ref{lem4.2}, $H$ is a $(0,b)$-regular set of $GC(G,S,\alpha)$.

(a)$\Rightarrow$(c).
From the construction in (a)$\Rightarrow$(b),  $\{H, H\ell_1,\dots,H\ell_{n-1}\}$ is a $(-b)$-equitable partition of $GC(G,S,\alpha)$, where each right coset   $H\ell_i$ is a $(0,b)$-regular set of $GC(G,S,\alpha)$.
By Lemma \ref{lem3}, $b$ divides $|S|$.
For each $0\leq a\leq \frac{|S|}{b}-1$, define $H_a=\bigcup_{i=0}^a H\ell_i$ with $H\ell_0=H$.
A direct counting argument shows that each vertex in $H_a$ has $ab$ neighbors in $H_a$, and each vertex in $G\setminus H_a$ has $(a+1)b$ neighbors in $H_a$.
Thus $H_a$ is an $(ab,(a+1)b)$-regular set of $GC(G,S,\alpha)$.

(c)$\Rightarrow$(a).
Taking $a=0$  yields that $H$ is a $(0,b)$-regular set of $GC(G,S,\alpha)$.
\end{proof}

\begin{remark}

In  \cite{LiaoLiuZhang2025}, the connectivity of generalized Cayley graphs was characterized.
The  generalized Cayley graph $GC(G,S,\alpha)$ is connected if and only if $G=\langle S\rangle$, $|G:\langle SS^{-1}\rangle|\leq 2$, and $\alpha(\langle SS^{-1}\rangle)=\langle SS^{-1}\rangle$.
In fact, if $GC(G,S,\alpha)$ is connected and $\langle SS^{-1}\rangle\neq G$, then $GC(G,S,\alpha)$ is bipartite with bipartition $\langle SS^{-1}\rangle$ and $\langle SS^{-1}\rangle s$ for any $s\in S$.
Thus, the  subgroup $\langle SS^{-1}\rangle$ is a $(0,|S|)$-regular set of $GC(G,S,\alpha)$.
\end{remark}

\medskip
\noindent\textbf{Data availability}
No data was used for the research described in the paper.

\medskip
\noindent\textbf{Declarations}

\noindent
{\bf Conflict of interest} The authors declare that they have no conflict of interest to this paper.



\end{sloppypar}


\frenchspacing
\begin{thebibliography}{99}

\bibitem{AhanjidehAkhlaghiRodrigues2026}
N. Ahanjideh, Z. Akhlaghi, B.G. Rodrigues. On perfect codes and total perfect codes in Cayley sum graph. J. Algebraic Combin. 2026, 63(4): No52.

\bibitem{BaileyCameron2019}
R.A. Bailey, P.J. Cameron, A.L. Gavrilyuk,  S.V. Goryainov.
Equitable partitions of Latin-square graphs. J. Combin. Des. 2019, 27(3): 142--160.

\bibitem{ChenZhou2016}
G. Chen, S.M. Zhou. Perfect codes in Cayley graphs. SIAM J. Discrete Math. 2016, 30(3): 1508--1525.



\bibitem{ChenLiZhang2025}
Z.G. Chen, J.J. Li, J.Y. Zhang. Subgroup perfect codes in lie type simple groups of rank one. Des. Codes Cryptogr. 2025, 93(10): 4275-4292.

\bibitem{H.K.P2017}
A. Hujdurovi\'{c}, K. Kutnar, P. Petecki, et al.  {On automorphisms and structural properties of generalized Cayley graphs}.  Filomat. 2017, 31(13): 4033--4040.

\bibitem{KhaefiAkhlaghiKhosravi2025}
Y. Khaefi, Z. Akhlaghi, B. Khosravi. On the subgroup regular sets in Cayley graphs. J. Algebraic Combin. 2025, 61(3): No38.

\bibitem{LiaoLiuLiTian2026}
Q.F. Liao, W.J. Liu, Y.M. Li, D.L. Tian. Subgroup (total) perfect codes of generalized Cayley graphs. Bull. Malays. Math. Sci. Soc. 2026, 49(1): No41.

\bibitem{LiaoLiuZhang2025}
Q.F. Liao, W.J. Liu, P.L. Zhang. On structural properties of generalized Cayley graphs. Filomat. 2025, 39(30): 10871--10882.

\bibitem{MaFengWang2020}
X.L. Ma, M. Feng, K.S. Wang. Subgroup perfect codes in Cayley sum graphs. Des. Codes Cryptogr. 2020, 88(7): 1447--1461.

\bibitem{M.L.Z}
X.L. Ma, L. Li, G. Zhong. Perfect codes in proper intersection power graphs of finite groups. Appl. Algebra Engrg. Comm. Comput. 2025, 36(3): 543-555.

\bibitem{Marusic1992}
D. Maru\v{s}i\v{c}, R. Scapellato, N. Zagaglia Salvi. {Generalized Cayley graphs}. Discrete Math. 1992, 102: 279--285.


\bibitem{SeiedaliAkhlaghiKhosravi2026}
F.S. Seiedali, Z. Akhlaghi, B. Khosravi. Generalized Cayley graphs and perfect code. Discrete Math. 2026, 349(2): No114805.

\bibitem{SeiedaliKhosraviAkhlaghi2025}
F.S. Seiedali, B. Khosravi, Z. Akhlaghi. Regular Sets in Cayley Sum Graphs. Bull. Malays. Math. Sci. Soc. 2025, 48(2): No35.

\bibitem{WangXu2024}
X.M. Wang, S.J. Xu. On subgroup regular sets in Cayley sum graphs. Discrete Math. 2024, 347(9): No114059.

\bibitem{WangXuZhou2024}
X.M. Wang, S.J. Xu, S.M. Zhou. On regular sets in Cayley graphs.
J. Algebraic Combin. 2024, 59: 735--759.


\bibitem{Zhang2024}
J.Y. Zhang. On subgroup perfect codes in Cayley sum graphs. Finite Fields Appl. 2024, 95: No102393.


\bibitem{ZhangYuan2026}
J.Y. Zhang, W. Yuan. On regular sets in Cayley sum graphs. J. Algebra Appl. 2026, 25(07): No2650062.

\bibitem{ZhuLiuYang2023}
X.-M. Zhu, W.J. Liu, X. Yang. The isomorphism of generalized Cayley graphs on finite non-abelian simple groups. Discrete Math.  2023, 346(4): No113292.
\end{thebibliography}
\end{document}